\documentclass[a4paper,leqno,12pt]{article}
\usepackage{amssymb,xspace}
\usepackage{amsmath}
\usepackage{amsthm}

\swapnumbers
\newtheorem{num}{}[section] 

\newcommand{\beq}{\addtocounter{num}{1}\begin{equation}}
\newcommand{\alku}{\begin{num}}
\newcommand{\loppu}{\end{num}}

\newcommand{\lau}{{\bf Theorem. }}
\newcommand{\lem}{{\bf Lemma. }}
\newcommand{\pf}{{\it Proof. }\rm }
\newcommand{\rn}{{\sf R}^n}
\newcommand{\sub}{\subset}
\newcommand{\id}{{\rm id}}
\renewcommand{\span}{{\rm span}\,}
\newcommand{\aff}{{\rm aff}\,}
\newcommand{\cd}{{\rm cd}\,}
\newcommand{\tyh}{\varnothing}
\newcommand{\sm}{\setminus}
\newcommand{\iy}{\infty}
\newcommand{\imp}{$\Rightarrow$\xspace}
\newcommand{\eni}{$\varepsilon$-near\-isometry\xspace}
\def\pav[#1[{\left[#1\right[}
\def\ppav]#1]{\left]#1\right]}

\begin{document}

\title{A survey of nearisometries}
\author{Jussi V\"ais\"al\"a}
\date{}
\maketitle


\small
This paper appeared in Papers on Analysis, a volume dedicated to Olli Martio
on the occasion of his 60th birthday,
Report. Univ. Jyv\"askyl\"a 83, 2001, 305--315.

An addendum was added on p. 12 on January 1, 2002.
\normalsize

\section{Introduction}

\alku
\label{1.1} Notation. \rm Throughout this article, $E$ and $F$ are real
Banach spaces (sometimes Hilbert spaces or just euclidean spaces) of
dimension at least one. The norm of a vector  $x$ is written as $|x|$. In a
Hilbert space we let $x \cdot y$ denote the inner product of vectors $x$ and
$y$.

Open and closed balls with center $x$ and radius $r$ are written as $B(x,r)$
and $\bar B(x,r)$, respectively, and we use the abbreviations $B(r) =
B(0,r)$ and $\bar B(r) = \bar B(0,r)$. The open unit ball in the euclidean
space $\rn$ is $B^n = B(1)$.

We shall consider maps $f\colon A \to F$ where $A \sub E$. Without further
notice we shall always assume that $A \ne\tyh$. If $g\colon A \to F$ is
another map, we set
\[ d(f,g) = \sup \{  |fx - gx|: x \in A\}
\] with the possibility $d(f,g)=\iy$. We frequently consider the case where
$g = T|A$ is the restriction of an iso\-metry $T\colon E \to F$. Then we
simply write $d(T,f) = d(T|A,f)$. To simplify expressions we often omit
parentheses writing $fx = f(x)$ etc.

\loppu

\alku
\label{1.2} Nearisometries. \rm Let $A \sub E$. A map $f\colon A \to F$ is a
\textit{nearisometry} if there is a number
$\varepsilon \ge 0$ such that
\beq
\label{1.3} |x-y| - \varepsilon \le |fx-fy| \le |x-y| + \varepsilon
\end{equation} for all $x,y \in A$. More precisely, we say that such a map
is an
$\varepsilon$-\textit{nearisometry}. Observe that $f$ need not be
continuous. In the literature, the nearisometries and
$\varepsilon$-nearisometries are often called approximate isometries and
$\varepsilon$-isometries, respectively.

The condition (\ref{1.3}) can be regarded as a perturbation of the
iso\-metry condition
$|fx-fy| = |x-y|$. Another kind of perturbation is given by the
multiplicative condition
\beq\label{1.4} |x-y|/M \le |fx-fy| \le M|x-y|.
\end{equation} A map satisfying (\ref{1.4}) is called
$M$-\textit{bilipschitz}.

Neither of the conditions (\ref{1.3}) and (\ref{1.4}) implies the other.
Condition (\ref{1.3}) is stronger than (\ref{1.4}) for large distances but
weaker for small distances. A bilipschitz map is always continuous, even an
embedding. 

However, if $A$ is bounded, then a $(1+t)$-bilipschitz map $f \colon A \to
F$ is a
$td(A)$-nearisometry, where $d$ denotes diameter. In the other direction, an
$\varepsilon$-nearisometry satisfies (\ref{1.4}) with $M = 1+\varepsilon$
whenever $|x-y| \ge 1+\varepsilon$.  The basic question considered in this
paper is the following stability problem:

 Given an
$\varepsilon$-nearisometry $f\colon A \to F$, does there exist an isometry
$T\colon E \to F$ such that $d(T,f)$ is small if $\varepsilon$ is small? In
other words, are the nearisometries near isometries? 

It turns out that in many cases we can even find a linear bound $d(T,f) \le
c\varepsilon$. This is true, for example, if
$A$ is the whole space and $f$ is surjective. For Hilbert spaces this was
proved in 1945 by D.H. Hyers and S.M. Ulam in the famous paper \cite{HU1}.

The classical case $A = E$ is considered in Section 2. In Section 3 we
consider the case where
$A$ is an unbounded subset of $E$, and Section 4 deals with bounded sets $A
\sub E$. In Section 5 we mention some related results.
\loppu

\section{Whole spaces}

\alku
\label{2.1} Preliminary results. \rm We recall the classical Mazur-Ulam
theorem (see \cite[14.1]{BL}): Every surjective isometry
$T\colon E \to F$ is affine, that is, $Tx = Sx + T(0)$ where $S$ is linear.
This is also true for nonsurjective isometries if $F$ is strictly convex, in
particular, if $F$ is Hilbert.

In this section we consider $\varepsilon$-nearisometries $f\colon E \to F$
defined in the whole space $E$. There is an extensive literature dealing
with this case starting with the influential paper \cite{HU1} of Hyers and
Ulam; see \ref{his}. To find an isometric approximation $T\colon E
\to F$ we may without loss assume that $f(0) = 0$. If $T$ is linear and if
$d(T,f)\ < \iy$, it is easy to see that $Tx = \lim_{t\to\iy} f(tx)/t$ for
each $x \in E$. Consequently, if $T$ is surjective or if $F$ is strictly
convex, an isometric approximation of $f$ is uniquely determined up to
translation.

We recall that the \textit{Jung constant} $J(E)$ of $E$ is the infimum of
all numbers $r>0$ such that every set $A \sub E$ of diameter $d(A) \le 2$
can be covered by a ball of radius $r$. We have always $1 \le J(E) \le 2$,
and $J(\rn) = \sqrt{2n/(n+1)} < \sqrt 2$ by the classical result proved by
H.W.E. Jung \cite{Ju} in 1901; see \cite[2.10.45]{Fe}. Moreover, $J(E) =
\sqrt 2$ for infinite-dimensional Hilbert spaces; see \cite[Th. 2]{Da} or
\cite[p. 704]{Se1}. The upper bound
$J(E) = 2$ is obtained, for example, by the Banach spaces $c$ and $c_0$.

\loppu

The following result summarizes the work of several authors during
1945--1998.

\alku
\label{fund} {\bf Fundamental theorem.} Suppose that $f\colon E \to F$ is a
surjective
$\varepsilon$-near\-isometry with $f(0)=0$. Then there is a surjective
linear iso\-metry
$T\colon E \to F$ with $d(T,f) \le 2\varepsilon$. The bound is the best
possible.

Moreover, for each $t>0$ there is a surjective (hence affine) iso\-metry $S\colon E \to F$
with $d(S,f) \le J(E)\varepsilon + t$. Also here the bound
$d(S,f) \le 2\varepsilon$ is the best possible.
\loppu

\alku
\label{his} History. \rm The first part of the Fundamental theorem was
proved for Hilbert spaces by Hyers and Ulam
\cite{HU1} in 1945 with the bound $10\varepsilon$. Surprisingly, it took 38
years until J. Gevirtz \cite{Ge} obtained a proof for all Banach spaces
(with the bound $5\varepsilon$) in 1983. Meanwhile, various special cases
were considered in \cite{HU2}, [Bo1--Bo5] and \cite{Gr}. P.M. Gruber
\cite{Gr} got quite close to the solution in 1978. For example, he proved
that if there is $T$ with $d(T,f) < \iy$, then $d(T,f) \le 5\varepsilon$.
Moreover, he obtained the result for all finite-dimensional Banach spaces.

M. Omladi\v c and P. \v Semrl \cite{OS} improved the ideas of Gruber and
obtained the bound
$2\varepsilon$ in 1995. By a simple example they also proved the sharpness
of the bound. The following example, due to S.J. Dilworth \cite{Di}, is
still simpler.

Let $\varepsilon > 0$ and define $f\colon {\sf R} \to {\sf R}$ by $fx =
x-\varepsilon$ for $x
\ne 0,\varepsilon$ and by $fx = -x$ for $x=0,\varepsilon$. This map is
obviously an \eni. If
$T\colon {\sf R}
\to {\sf R}$ is a linear isometry, then either $T = \id$ or $T = -\id$.
Hence $d(T,f)$ is
$2\varepsilon$ or
$\iy$. 

The second part of the theorem follows rather easily from the first part, as
shown by \v Semrl
\cite{Se1} in 1998. In the same paper he also proved the sharpness of the
bound $2\varepsilon$ in the second part for the space $E = c$. A somewhat
simpler proof for the space $c_0$ is given in \cite{HV}. We remark that one 
can choose $t=0$ if each set $A \sub E$ with $d(A) \le 2$ can be covered by 
a ball of radius $J(E)$. This is true, for example, if $E$ is a Hilbert
space.

It is natural to ask whether the bound $d(S,f) \le J(E)\varepsilon + t$ is
sharp for every Banach space $E$. An affirmative answer for Hilbert spaces
was given in \cite{HV} but the general case is open.

One might also think that the result holds with better bounds if only
continuous maps $f\colon E
\to F$ are considered. However, this is not the case at least for Hilbert
spaces, as proved in
\cite{HV}. 

\loppu

\alku
\label{2.4} Remarks on the proof. \rm The proof of the Fundamental theorem
\ref{fund} is entirely elementary, and it can be understood (at least if $f$
is bijective) by anyone who knows the definition and some basic properties of
Banach spaces. A nice presentation is given in the book \cite[15.2]{BL} of
Y. Benyamini and J. Lindenstrauss. The proof consists of three steps. 

1.  Given $x \in E$, we show that the sequence of the points $y_n =
2^{-n}f(2^nx)$ is Cauchy and hence converges to a limit $Tx$. We thus obtain
a map $T\colon E \to F$, which is clearly an isometry.

2. We show that $T$ is surjective.

3. The inequality $|Tx-fx| \le 2\varepsilon$ is proved.

If the spaces $E$ and $F$ are Hilbert, the proof is substantially easier
than in the general case. For example, the sequence $(y_n)$ converges
whenever $f$ is a near\-isometry defined in the set $\{  2^nx: x \in {\sf
N}\}$. In particular, the surjectivity of $f$ is not needed in Step 1. A
simplified proof for the Hilbert space case is given in \cite[5.1]{unbdd}

However, the surjectivity condition of \ref{fund} cannot be omitted even in
the case $f\colon {\sf R} \to {\sf R}^2$; see \ref{2.5}. It can be omitted
if $E$ and $F$ have the same finite dimension (Theorem \ref{2.6}). This was
proved by R. Bhatia and  \v Semrl \cite{BS} for euclidean spaces and by
Dilworth \cite{Di} in the general case.
\loppu

\alku
\label{2.5} Example. \rm Let $\varepsilon>0$ and define $f\colon {\sf R} \to
{\sf R}^2$ by $fx = (x, \sqrt{2\varepsilon |x|})$. It is easy to see that
$f$ is an $\varepsilon$-near\-isometry. However $d(T,f) = \iy$ for every
isometry $T\colon {\sf R} \to {\sf R}^2$.
\loppu

\alku
\label{2.6}
\lau
 If $\dim E = \dim F < \iy$, then the Fundamental theorem holds without the
surjectivity condition.
\loppu

\alku
\label{2.7} Relaxing the surjectivity condition. \rm Although the
surjectivity condition of the Fundamental theorem \ref{fund} cannot be
removed, it can be weakened. For example, it can be replaced by the
condition that $F \sm fE$ is bounded. A more general result is given in
\ref{cobounded}.

A still weaker condition is given by the property of
$\delta$-\textit{ontoness}. This can also be considered relative to a closed
linear subspace $F_1$ of $F$. Let $\delta \ge 0$.  We say that a map
$f\colon E \to F$ is $\delta	$-\textit{onto} $F_1$ if the Hausdorff distance
$d_H(fE,F_1)$ is at most $\delta$, that is, $d(fx,F_1) \le \delta$ and
$d(y,fE) \le \delta$ for all $x \in E$ and
$y \in F_1$. In particular, a map $f$ is $\delta$-onto $F$ iff $F \sm fE$
contains no ball of radius larger than $\delta$.
 \loppu

\alku
\label{2.8}
\lau Suppose that $f\colon E \to F$ is an $\varepsilon$-near\-isometry with
$f(0)=0$ and that $f$ is $\delta$-onto $F_1$ where $F_1$ is a closed linear
subspace of $F$. Then there is a surjective linear iso\-metry $T\colon E \to
F_1$ such that $d(T,f) \le 2\varepsilon + 4\delta$.
\loppu

This result was obtained by Dilworth \cite{Di} with a somewhat larger bound.
The bound
$2\varepsilon+4\delta$ is from \cite{SV.prep}. See also \cite[Th. 3]{Ta}.
The first term
$2\varepsilon$ is the best possible by the Fundamental theorem but the
second term $4\delta$ is presumably not. It must be at least $\delta$ as is
seen from the following example; see
\cite[p. 473]{Di}:

 Let $\delta > 0,\ M > 0$, and let $g\colon {\sf R} \to {\sf R}$ be an
increasing $M$-Lipschitz function with $\lim_{x\to -\iy} gx = 0,\ \lim_{x\to
\iy} gx = 1$. Define $f\colon {\sf R} \to {\sf R}^2$ by $fx = (x,\delta
gx)$.  Then
$f$ is an \eni with $\varepsilon = M\delta/2$. Indeed, if $a \in {\sf R},\ b
= a+d > a,\ d' = |fb - fa|,\ t = \delta(gb-ga)$, then ${d'}^2 = d^2 + t^2$
and
\[ 0 \le d' - d = t^2/(d + d') \le Md\delta/2d = M\delta/2.
\] The map
 $f$ is $\delta$-onto $F_1 = {\sf R} \times \{  0\}$, and $d(T,f) = \delta$
and $d(T,f) = \iy$ for the two linear isometries $T\colon  {\sf R} \to F_1$.
Since $\varepsilon$ is arbitrarily small, Theorem \ref{2.8} is not true if
$4\delta$ is replaced by a number less than
$\delta$.

\alku
\label{2.9} The case $F_1 = F$. \rm This case of \ref{2.8} is substantially
different from the case $F_1 \ne F$. One can show that the bound
$2\varepsilon + 4\delta$ can then be replaced by $2\varepsilon + 2\delta$
but the following conjecture looks plausible:

\loppu

\alku
\label{2.10} Conjecture. \rm The Fundamental theorem is true if the
surjectivity of $f\colon E \to F$ is replaced by the condition that $f$ is
$\delta$-onto $F$ for some $\delta$.
\loppu

\alku
\label{2.11} Results. \rm ({\bf See Addendum at the end of the paper.)} At the
time of this writing (August 2001), the conjecture is open. It is known to be
true in the following cases:

(1) The norm of $E$ is Fr\' echet differentiable in a dense set \cite[Th.
2]{Di}. This class includes all Asplund spaces (separable subspaces have
separable duals) and hence all reflexive spaces.

(2) $E = F $ is $l_p$ or $ L_p(X,\mu)$, where $1 \le p \le \iy$ and
$(X,\mu)$ is a measure space
\cite{SV.prep}. The case $1 < p < \iy$ is included in (1).

(3) $E = F = C(X)$, the space of continuous real-valued functions in a
compact space $X$
\cite{SV.prep}. 

Moreover, the conjecture is true for all spaces with the bound
$2\varepsilon$ replaced by
$3\varepsilon$ \cite{SV.prep}.
\loppu

\alku
\label{2.12} The number $\tau(Q)$. \rm We next try to replace the
$\delta$-ontoness in \ref{2.10} by a still weaker condition. Let $Q$ be a
nonempty set in a Banach space and let $u$ be a unit vector. We set
\[
\varrho(u,Q) = \liminf_{|t|\to \iy} d(tu,Q)/|t|,\quad \tau(Q) = \sup_{|u|=1}
\varrho(u,Q).
\] Then
\[
\varrho(-u,Q) = \varrho(u,Q),\quad 0 \le \varrho(u,Q) \le 1, \quad 0 \le
\tau(Q) \le 1.
\] Moreover, these numbers are translation invariants: $\varrho(u,Q + z) =
\varrho(u,Q),\ \tau(Q+z) = \tau(Q)$.

If $Q$ is bounded or contained in a hyperplane, then $\tau(Q) = 1$. If $Q$
contains a half space, then $\tau(Q) = 0$. If the function
$y
\mapsto d(y,Q)$ is bounded, then $\tau(Q) = 0$. Hence $\tau(fE) = 0$ if a
map $f\colon E \to F$ is
$\delta$-onto $F$ for some $\delta$. 

In view of Conjecture \ref{2.10}, it is reasonable to ask whether an
$\varepsilon$-near\-isometry $f\colon E \to F$ with $\tau(fE)= 0$ or maybe
with $\tau(fE) < 1$ can be approximated by a surjective iso\-metry. An
answer for Hilbert spaces is given in
\ref{tauhilbert}. Observe that $\tau(f{\sf R}) = 1$ for the map $f\colon
{\sf R}\to {\sf R}^2$ considered in \ref{2.5}.

\loppu

\alku
\label{tauhilbert}
\lau {\rm \cite[5.4]{unbdd}} Suppose that $E$ and $F$ are Hilbert spaces and
that $f\colon E \to F$ is an
$\varepsilon$-near\-isometry with $f(0)=0$ and $\tau(fE)<1$. Then there is a
surjective linear iso\-metry $T\colon E \to F$ with $d(T,f) \le
2\varepsilon$.
\loppu

Observe that the condition $d(T,f) \le 2\varepsilon$ implies that $\tau(fE)
= 0$. It follows that $\tau(fE) \in \{  0,1\}$ for each near\-isometry
between Hilbert spaces.

\alku
\label{2.15}  Inverse results. \rm
 An isometry $f\colon E \to F$ between Banach spaces with $f(0)=0$ need not
be linear, but T. Figiel (\cite{Fi}, \cite[14.2]{BL}) and W. Holszty\' nski
\cite{Ho} proved in 1968 that there is a unique linear map $T\colon 
\overline{{\rm span}}\, fE \to E$ such that $|T| = 1$ and $T  f = \id$. In
view of this result, it is reasonable to conjecture that if $f\colon E \to
F$ is an \eni with $f(0)=0$, then there is a linear map $T\colon 
\overline{{\rm span}}\, fE \to E$ such that $|T|=1$ and $d(T f,\id) \le
c\varepsilon$ for some constant $c$. This conjecture was disproved by S.
Qian \cite[Ex. 1]{Qi} in 1995, but Qian proved that the conjecture holds
(with $c=6)$ in certain cases, for example, if $F$ is Hilbert or if $E$ and
$F$ are $L^p$ spaces with $1 < p < \iy$. Further results in this direction
are given in \cite{SV.prep}.
\loppu

\section{Unbounded subsets}

In this section we consider $\varepsilon$-near\-iso\-metries $f\colon A \to
F$ where $A$ is an unbounded subset of $E$. We look for a surjective
iso\-metry $T\colon E \to F$ such that $d(T,f)$ is finite and hopefully
bounded by $c\varepsilon$ for some constant $c$, possibly depending on
$A$. As in Section 2, the proofs are based on the behavior of $f$ near the
point at infinity.

The results in this section are due to the author \cite{unbdd}. Most of them
deal with Hilbert spaces or just with euclidean spaces. However, the
following improvement of the Fundamental theorem \ref{fund} is valid for all
Banach spaces. It shows that the Fundamental theorem is not, after all, a
truly global result but a local property of maps at the point
$\iy$. The result was suggested to the author by O. Martio.

\alku
\label{cobounded}
\lau Suppose that $A \sub E$ and that $f\colon A \to F$ is an
$\varepsilon$-near\-isometry such that the sets $E \sm A$ and $F \sm fA$ are
bounded. Then there is a surjective iso\-metry
$T\colon E \to F$ with $d(T,f) \le 2\varepsilon$. For each $x_0\in A$ we can
choose $T$ so that
$Tx_0 = fx_0$.
\loppu

\pf We may assume that $x_0 = 0,\ fx_0 = 0$. Choose a number $R >
\varepsilon$ such that $E \sm A
\sub B(R)$ and $F \sm fA \sub B(2R)$. Define $f_1\colon E \to F$ by $f_1x =
2x$ for $|x| \le R$ and by
$f_1x = fx$ for $|x| > R$. Then $f_1$ is a surjective near\-isometry. By the
Fundamental theorem
\ref{fund}, there is a linear surjective isometry
$T\colon E \to F$ with $d(T,f_1) < \iy$. Since $d(T,f) < \iy$, an easy
modification of the proof of
\ref{fund} (see
\cite[p. 362]{BL}) shows that $d(T,f) \le 2\varepsilon$.  $\square$ 

\alku
\label{3.2} Hilbert spaces. \rm In the rest of this section we assume that
$E$ and $F$ are Hilbert spaces. As mentioned in 
\ref{2.4}, the sequence of points $y_n = 2^{-n} f(2^nx)$ converges as soon
as the near\-isometry
$f$ is defined in the set $\{  2^nx: n \in {\sf N}\}$. A more general result
is given in
\ref{3.3}. These results are not valid in general Banach spaces, which makes
the theory in Banach spaces considerably more difficult.
\loppu

\alku
\label{3.3}
\lem Suppose that $E$ and $F$ are Hilbert spaces and that $(x_j)$ is a
sequence in $E$ such that
$|x_j| \to\iy$ and $x_j/|x_j|$ converges to a limit as $j \to\iy$. Suppose
also that $f\colon
\{  x_j:j \in {\sf N}\} \to F$ is a near\-isometry. Then the sequence of the
points $fx_j/|x_j|$ is convergent.
\loppu

\alku
\label{3.4} Terminology. \rm Suppose that $A$ is an unbounded subset of a
Hilbert space $E$. We say that a unit vector $u \in E$ is a \textit{cluster
direction} of $A$ if there is a sequence $(x_j)$
 in $A$ such that $|x_j| \to\iy$ and $x_j/|x_j| \to u$. We let $\cd A$
denote the set of all cluster directions of $A$. 

The theory is easier in the case $\dim E < \iy$, because then each sequence
$(x_j)$ with $|x_j|
\to\iy$ has a subsequence $(y_j)$ such that the sequence $(y_j/|y_j|)$ is
convergent. In particular, $\cd A$ is then compact and nonempty. In a
general Hilbert space, an unbounded set need not have any cluster directions.

For an unbounded set $A \sub E$ we define a number $\mu(A) \in [0,1]$ as
follows. First, for each unit vector $e \in E$ we set $\sigma(e,A) = \sup
\{  |u \cdot e| : u \in \cd A\}$ and then
\[
\mu(A) = \inf \{  \sigma(e,A): |e|=1\}.
\] The number $\mu(A)$ is small iff $\cd A$ lies in a narrow neighborhood of
a hyperplane through the origin.

For example, $\mu(A) = 1$ if $A$ contains a half space or if $d(x,A)$ is
bounded for $x \in E$.
 Furthermore, $\mu(C) = \sin \alpha$ for the cone $C = \{  x: |x
\cdot  e| \ge |x| \cos \alpha\}$ where
$0 <\alpha \le \pi/2,\ |e|=1$.

If $\dim E < \iy$, then $\mu(a)$ is the infimum of all
$t \in [0,1]$ such that the double cone $\{  x \in E: |x \cdot e| > t|x| \}$
meets $A$ in a bounded set for some unit vector $e$. Moreover, in this case
$\mu(A) = \sqrt{1-
\tau(A)^2}$, where $\tau$ is defined in \ref{2.12}.

Let $c > 0$. We say that a set $A \sub\rn$ has the $c$-\textit{isometric
approximation property}, abbreviated $c$-IAP, if for every $\varepsilon \ge
0$ and for every \eni $f\colon A
\to\rn$ there is an isometry $T\colon \rn \to \rn$ with $d(T,f) \le
c\varepsilon$. If $A$ has the $c$-IAP and contains at least two points, then
clearly $c \ge 1/2$.

For example, the whole space $\rn$ has the $\sqrt 2$-IAP by \ref{2.6}. From
Example \ref{2.5} it follows that a line in $\rn,\ n \ge 2$, does not have
the $c$-IAP for any $c$. The following result
\cite[2.3]{unbdd} gives a quantitative geometric characterization for
unbounded subsets of
$\rn$ with the $c$-IAP. The corresponding question for bounded sets will be
considered in  Section 4.

\loppu

\alku
\label{3.5}
\lau For an unbounded set $A \sub \rn$, the following conditions are
quantitatively equivalent. 

$(1)$ $A$ has the $c$-{\rm IAP},

$(2)$ $\mu(A) \ge 1/c' > 0$.

\noindent More precisely, $(1)$ implies $(2)$ with $c' = 17c$, and $(2)$
implies $(1)$ with $c =
\sqrt 2 c'$. The constant $\sqrt 2$ is the best possible.
\loppu

\alku
\label{3.6} Remarks. \rm
 1. It follows from \ref{3.5} that $\mu(A) = 0$ iff $A$ does not have the
$c$-IAP for any
$c$. This happens, for example, if $A$ is a linear subspace of $\rn$ with
$\dim A < n$. This can be directly seen as in Example \ref{2.5}. In fact,
the proof for the part (1) \imp (2) of
\ref{3.5} is based on an elaboration of this example.

2. We sketch the proof for the part (2) \imp (1). Suppose that $\mu(A) \ge
1/c'$ and that
$f\colon A \to \rn$ is an \eni. We may assume that $0 \in A$ and that
$f(0)=0$. We first define a map $\varphi\colon \cd A \to \rn$ as follows.
Let $u \in \cd A$ and choose a sequence $(x_j)$ in $A$ such that $|x_j|
\to\iy$ and $x_j/|x_j| \to u$. By \ref{3.3}, the limit $\varphi u =
\lim_{j\to\iy} fx_j/|x_j|$ exists. One can show that the limit is
independent of the choice of the sequence and that the map $\varphi\colon
\cd A \to \rn$ is an iso\-metry, which extends to a linear iso\-metry $T
\colon \rn\to\rn$ with $d(T,f) \le 2\varepsilon$. The bound $\sqrt 2
\varepsilon$ is obtained by composing $T$ with a suitable translation.

Most of these arguments can be carried out in an arbitrary Hilbert space,
and we obtain the following variation of the part (2) \imp (1) of \ref{3.5}:

\loppu

\alku
\label{3.7}
\lau Suppose that $E$ and $F$ are Hilbert spaces and that $A \sub E$ is an
unbounded set with $\mu(A)
\ge 1/c'$. Suppose also that $f\colon A \to F$ is an \eni. Then there is an
iso\-metry $T\colon E
\to F$ such that $d(T,Pf) \le \sqrt 2 \varepsilon$, where $P\colon F \to TE$
is the orthogonal projection onto the affine subspace $TE$ of $F$.
\loppu

\section{Bounded subsets}

In this section we consider near\-iso\-metries of a bounded set $A \sub\rn$.
The target space will be $\rn$ (with the same $n$) except in the first
result \ref{4.1}. Most of the results are due to P. Alestalo, D.A. Trotsenko
and the author.

The proofs are very much different from those in Sections 2 and 3, because
we cannot use a limiting process where points tend to $\iy$. The basic tool
is the simple formula
\[ 2a \cdot b = |a|^2 + |b|^2 - |a-b|^2,
\] and the proofs are elementary but rather long except if the set $A$ is
sufficiently regular (for example a ball), in which case one can make use of
the ideas of F. John \cite{Jo}; see \ref{4.3}.

We identify $\rn$ in the natural way with a subspace of the Hilbert space
$l_2$.

\alku
\label{4.1}
\lau {\rm \cite[2.2]{ATV1}} Suppose that $A \sub\rn$ is bounded and that
$f\colon A \to l_2$ is an
$\varepsilon d(A)$-near\-isometry with $\varepsilon \le 1$. Then there is a
surjective isometry $T\colon l_2 \to l_2$ with $d(T,f) \le c(n)\sqrt
{\varepsilon} d(A)$. If $fA \sub\rn$, we can choose $T$ so that $T\rn \sub
\rn$.
\loppu

\alku
\label{4.2} Example. \rm Let $A$ be the interval $[-1,1]$, let $\varepsilon
> 0$, and let $f\colon A \to {\sf R}^2$ be the map
$fx = (x,|x|\sqrt\varepsilon)$. Then $f$ is an \eni, and $d(T,f) \ge \sqrt
\varepsilon/2$ for each iso\-metry $T\colon A \to {\sf R}^2$. Hence the
bound in \ref{4.1} has the correct order of magnitude.
\loppu

\alku
\label{4.3} John's method. \rm  F. John \cite{Jo} considered in 1961
isometric approximation of locally
$(1+\varepsilon)$-bilipschitz maps $f\colon G \to\rn$ where $G \sub \rn$ is
a ball or, more generally, of a class later called John domains. His method
is elegant compared with the proofs of the other results in this section. It
can easily be modified so as to prove the IAP (see
\ref{3.4}) of sufficiently regular bounded sets. We give the result for
balls; a more general result is given in \cite[3.12]{ATV1}.
\loppu

\alku
\label{4.4}
\lau A ball in $\rn$ has the $c$-{\rm IAP} with $c = 10n^{3/2}$.
\loppu

\alku
\label{4.5} Remark. \rm The constant $c$ in \ref{4.4} must depend on $n$.
This was recently proved by E. Matou\v skov\' a \cite{Ma}, who showed that
for each $t > 0$ there is an integer $n$ and a
$(1+t)$-bilipschitz map $f\colon \rn\to\rn$ such that $d(T|B^n,f|B^n) \ge
1/\sqrt 2$ for each iso\-metry $T\colon \rn\to\rn$. Since $f|B^n$ is a
$2t$-near\-isometry, the unit ball $B^n$ does not have the $c$-IAP for $c <
(2t\sqrt 2)^{-1}$.

The map $f$ is defined as follows. Let $h\colon {\sf R}^2 \to {\sf R}^2$ be
the
$(1+t)$-bilipschitz spiral map defined by $h(r,\varphi) = (r, \varphi +
\pi/2 + t\log r)$ in polar coordinates. Choose an integer $N > \pi/t\log 2$
and set
$m = 2^N,\ n = 2m$. Then $\rn = E_1 \oplus \dots \oplus E_m$ where $E_j =
\span (e_j,e_{m+j})$. The map $h$ induces in a natural way maps $h_j\colon
E_j \to E_j$, and we set $h = h_1
\oplus\dots\oplus h_m\colon \rn\to\rn$. Then $f$ is $(1+t)$-bilipschitz with
$|fx| = |x|$ and $f(-x) = fx$ for all $x \in \rn$. One can show that the
image of ${\sf R}^m = \span (e_1,\dots,e_m)$ contains an orthonormal basis
$\bar u$ of $\rn$. If $T\colon \rn\to\rn$ is an iso\-metry, then
$L = T{\sf R}^m$ is an affine subspace of $\rn$, and one can show that there
is a member $u$ of
$\bar u$ such that $d(u,L) \vee d(-u,L) \ge 1/\sqrt 2$. Hence $f$ is the
desired map.
\loppu

\alku
\label{4.6} Thickness. \rm For a unit vector $u \in\rn$ we define the
projection $\pi_u\colon \rn \to {\sf R}$ by $\pi_u x = x \cdot u$. The
\textit{thickness} of a bounded set $A \sub\rn$ is the number
\[
\theta(A) = \inf_{|u|=1} d(\pi_uA).
\] We have always $\theta(A) \le d(A)$, and $\theta(A) = 0$ if and only if
$A$ is contained in a hyperplane. 

It follows from Example \ref{4.2} that a line segment $J \sub {\sf R}^2$
does not have the IAP. In this case we have $\theta(A) = 0$. One can show
that a bounded set $A \sub\rn$ containing at least $n+1$ points  has the IAP
if and only if $\theta(A) > 0$. If $A$ contains no isolated points, this
holds in the following quantitative form:
\loppu

\alku
\label{4.7}
\lau Suppose that $A \sub\rn$ is a bounded set without isolated points. Then
the following properties are quantitatively equivalent:

$(1)$ $A $ has the $c$-{\rm IAP},

$(2)$ $\theta(A) \ge d(A)/c'$.

\noindent More precisely, $(1)$ implies $(2)$ with a constant $c' = c'(c,n)$
and vice versa.
\loppu

Part (2) \imp (1) was proved in \cite[3.3]{ATV1}, and it is true for all
bounded sets. Part (1)
\imp (2) follows from \ref{solar} below.

\alku
\label{4.8} Sets with isolated points. \rm If $A \sub\rn$ is a bounded sets
containing isolated points, the part (1) \imp (2) of \ref{4.7} does not hold
quantitatively. This is seen from the following example due to  Trotsenko.

Let $0 < t \le 1$ and let $A \sub{\sf R}^2$ be the three-point set $\{ 
0,e_1,te_2\}$. Then
$d(A)/\theta(A) \ge t + 1/t$ is arbitrarily large. However, a direct proof
shows that $A$ has the 8-IAP for all $t$. 

As another example we consider the set $A' = \{  0,e_1,te_2,e_1+te_2\} \sub
{\sf R}^2$. It turns out that $A'$ has the $c$-IAP with some $c = c(t)$ but
$c(t)\to\iy$ as $t \to 0$. This is seen by considering the map $f\colon A'
\to {\sf R}^2$ with $f(e_1+te_2) = e_1-te_2$ and $fx = x $ for the other
three points $x \in A'$.

To get a quantitative geometric characterization for all bounded sets with
the $c$-IAP we introduce the following concept.
\loppu

\alku
\label{4.9} Definition. \rm Let $c \ge 1$. We say that a bounded set $A
\sub\rn$ is a $c$-\textit{solar system} if there is a finite set $H = \{ 
u_0,\dots,u_n\} \sub A$ such that

(1) $|u_k - u_0| \le c d(u_k, \aff (H \sm \{ u_k\})$ for all $1 \le k \le n$,

(2) $A \sm H \sub \bar B(u_0, c \min \{  |u_k - u_0|: 1 \le k \le n\})$.

\noindent Here $\aff S$ denotes the affine subspace spanned by a set $S
\sub\rn$.

The point $u_0$ plays a special role, and it is called the center of the
system. The points
$u_1,\dots,u_n$ are the planets and the set $A \sm H$ is the sun. The system
may degenerate to
$\{  u_0,\dots,u_i\}$ with $0 \le i < n$; then we assume that $u_{i+1} =
\dots = u_n = u_0$. If
$n \ge 3$ and if the system is nondegenerate, the planets do not lie in a
plane (as in the real solar system). Observe that there are no restrictions
for the distances $|u_k - u_0|$.

The three-point set $A$ of \ref{4.9} is a 1-solar system but the set $A'$ is
a $c$-solar system only for $c \ge \sqrt{1 + 1/t^2}$.
\loppu

\alku
\label{solar}
\lau {\rm \cite[2.5]{solar}} For a bounded set $A \sub\rn$, the following
conditions are quantitatively equivalent:

$(1)$ $A$ has the $c$-{\rm IAP},

$(2)$ $A$ is a $c'$-solar system.
\loppu

\section{Related results}

\alku
\label{weak} Weak near\-isometries. \rm Let $\varphi\colon
\pav[0,\iy[\to{\sf R}$ be an increasing function. We consider maps
$f\colon E \to F$ between Banach spaces satisfying the condition
\[
 \big| |fx-fy| - |x-y| \big| \le \varphi(|x-y|)
\] for all $x,y \in E$. If $\varphi$ is the constant function $\varphi(t) =
\varepsilon$, this means that $f$ is an \eni. J. Lindenstrauss and A.
Szankowski (\cite{LS}, \cite[15.4]{BL}) proved that if $f(0)=0$ and if
$\varphi$ increases so slowly that 
\[
\int_1^\iy \frac{\varphi(t)}{t^2} dt < \iy,
\] then there is a surjective linear iso\-metry $T\colon E \to F$ such that
$|Tx-fx|/|x[ \to 0$ as $|x| \to \iy$.
\loppu

\alku
\label{5.1} Stability. \rm  The theory considered in this article is an
example of stability (see \cite[p. 63]{Ul}). We consider a class
$C$ of maps (isometries) $f\colon X\to Y$. Then we relax the definition and
get a larger class
$C^*$ of maps ($\varepsilon$-near\-iso\-metries) involving a parameter
$\varepsilon$. Then we ask how well we can approximate a member $f$ of $C^*$
by members $T$ of $C$. Instead of estimating the distance $d(T,f)$ it is
sometimes more convenient to consider maps $T\colon Y
\to X$ and the distance $d(T f, \id)$. 

Various stability theories are considered in the survey articles of D.H.
Hyers \cite{Hy} and G.L. Forti \cite{Fo}. We mention some examples.

1. $C$ = similarities, $C^*$ = quasisymmetric maps \cite[4.6]{ATV1}.

2.  $C$ = M\"obius maps, $C^*$ = quasiregular maps \cite[II.12.5]{Re}.

3. $C$ = additive maps, $C^*$ = almost additive maps \cite[15.1]{BL}.

4. $C$ = convex functions, $C^*$ = almost convex functions \cite{HU.convex}.
\loppu

\alku
\label{5.2} Applications. \rm The fundamental theorem \ref{fund} is
beautiful, but the author does not know of any applications of this result.
The IAP of thick sets (Th. \ref{4.7}) can be applied to prove the following
result on bilipschitz extensions \cite{ATV2}. Its proof follows the ideas in
\cite{blep}	 and
\cite{Tr}. 
\loppu

\alku
\label{5.3}
\lau For each positive integer $n$ and for each $c \ge 1$ there are positive
numbers $\varepsilon_0 =
\varepsilon_0(c,n)$ and $c' = c'(c,n)$ such that the following holds.

Suppose that $A$ is a subset of $\rn$ such that $\theta(A \cap B(x,r)) \ge
r/c$ whenever $x\in A$ and $A \sm B(x,r) \ne\tyh$. Then every
$(1+\varepsilon)$-bilipschitz map $f\colon A \to F$ with
$\varepsilon \le \varepsilon_0$ can be extended to a
$(1+c'\varepsilon)$-bilipschitz map
$g\colon \rn\to\rn$.
\loppu

\bigskip

{\bf Addendum.} January 1, 2002, after this paper appeared in Report. Univ.
Jyv\"askyl\"a 83.

 The paper
\cite{SV.prep}, which was in preparation at the time of writing this survey,
contains improvements to the subsections \ref{2.7}--\ref{2.11}. In particular,
Conjecture \ref{2.10} is true for all Banach spaces. Moreover, the estimate
$2\varepsilon + 4\delta$ of \ref{2.8} can be replaced by $2\varepsilon +
2\delta$ for all Banach spaces and by $2\varepsilon + \delta$ for Hilbert
spaces. Both bounds are sharp.

\noindent Preprints of the author can be downloaded from\\
www.helsinki.fi/$^{\sim}$jvaisala/preprints.html.
\bigskip

\noindent Matematiikan laitos\\ Helsingin yliopisto\\ PL 4, Yliopistonkatu
5\\ 00014 Helsinki, Finland\\
\texttt{jvaisala@cc.helsinki.fi}


\small

\begin{thebibliography}{xxx}

\bibitem [ATV1]{ATV1} P. Alestalo, D.A. Trotsenko and J. V\"ais\"al\"a,
Isometric approximation, Israel. J. Math. 125, 2001, 61--82.

\bibitem [ATV2]{ATV2} P. Alestalo, D.A. Trotsenko and J. V\"ais\"al\"a,
in preparation.

\bibitem [BL]{BL} Y. Benyamini and J. Lindenstrauss, Geometric nonlinear
functional analysis I, AMS Colloquium Publications 48, 2000.

\bibitem [B\v S]{BS}  R. Bhatia and P. \v Semrl, Approximate isometries on
Euclidean spaces, Amer. Math. Monthly 104, 1997, 497--504.

\bibitem [Bo1]{Bo1} D.G. Bourgin, Approximate isometries, Bull.  Amer. 
Math.  Soc. 52, 1946, 704--714.

\bibitem [Bo2]{Bo2} D.G. Bourgin, Approximately isometric and multiplicative
transformations on continuous function rings, Duke Math. J. 16, 1949,
385--397.

\bibitem [Bo3]{Bo3} D.G. Bourgin, Classes of transformations and bordering
transformations, Bull.  Amer.  Math.  Soc. 57, 1951, 223--237.

\bibitem [Bo4]{Bo4} R.D. Bourgin, Approximate isometries on finite
dimensional Banach spaces, Trans. Amer. Math. Soc. 207, 1975, 309--328.

\bibitem [Bo5]{Bo5} R.D. Bourgin, Two-dimensional $\varepsilon $-isometries. 
Trans.  Amer.  Math.  Soc. 244 (1978), 85--102. 

\bibitem [Da]{Da} J. Dane\v s, On the radius of a set in a Hilbert space,
Comment. Math. Univ. Carolin. 25, 1984, 355--362.

\bibitem [Di]{Di} S.J. Dilworth, Approximate isometries on
finite-dimensional normed spaces, Bull. London Math. Soc. 31, 1999, 471--476.

\bibitem [Fe]{Fe} H. Federer, Geometric measure theory, Springer, 1969.

\bibitem [Fi]{Fi} J. Fickett, Approximate isometries on bounded sets with an
application to measure theory, Studia Math. 72, 1982, 37--46.

\bibitem [Fig]{Fig} T. Figiel, On non linear isometric embeddings of normed
linear spaces, Bull. Acad. Polon. Sci. S\'er. Sci. Math. Astronom. Phys. 16,
1968, 185--188.

\bibitem [Fo]{Fo} G.L. Forti, Hyers-Ulam stability of functional equations
in several variables, Aequationes Math. 50, 1995, 143--190.

\bibitem [Ge]{Ge} J. Gevirtz, Stability of isometries on Banach spaces,
Proc. Amer. Math. Soc. 89, 1983, 633--636.

\bibitem [Gr]{Gr} P.M. Gruber, Stability of isometries, Trans. Amer. Math.
Soc. 245, 1978, 263--277.

\bibitem [Ho]{Ho} W. Holszty\' nski, Linearization of isometric embeddings
of Banach spaces. Metric envelopes, Bull. Acad. Polon. Sci. S\'er. Sci.
Math. Astronom. Phys. 16, 1968, 189--193.

\bibitem [HV]{HV} T. Huuskonen and J. V\"ais\"al\"a, Hyers-Ulam constants of
Hilbert spaces, preprint.



\bibitem [Hy]{Hy} D.H. Hyers, The stability of homomorphisms and related
topics, Global analysis - analysis on manifolds, ed. by T.M. Rassias,
Teubner, 1983, 140--150.
      
\bibitem [HU1]{HU1} D.H. Hyers and S.M. Ulam, On approximate isometries,
Bull. Amer. Math. Soc. 51, 1945, 288--292.

\bibitem [HU2]{HU2} D.H. Hyers and S.M. Ulam, Approximate isometries of the
space of continuous functions, Math. Ann. 48, 1947, 285--289.

\bibitem [HU3]{HU.convex} D.H. Hyers and S.M. Ulam, Approximately convex
functions, Proc. Amer. Math. Soc. 3, 1952, 821--828.

\bibitem [Jo]{Jo} F. John, Rotation and strain, Comm. Pure Appl. Math. 14,
1961, 391--413.

\bibitem [Ju]{Ju} H.W.E. Jung, \"Uber die kleinste Kugel, die eine
r\"aumliche Figur einschliesst, J. Reine Angew. Math. 123, 1901, 241--257.

\bibitem [LS]{LS}   J. Lindenstrauss and A. Szankowski, Non linear
perturbations of iso\-metries, Ast\'erisque 131, 1985, 357--371.

\bibitem [Ma]{Ma} E. Matou\v skov\'a, Almost isometries of balls, J. Funct.
Anal, to appear.

\bibitem [O\v S]{OS} M. Omladi\v c and P. \v Semrl, On nonlinear
perturbations of iso\-metries, Math. Ann. 303, 1995, 617--628.

\bibitem [Re]{Re} Yu. Reshetnyak, Space mappings with bounded distortion,
Amer. Math. Soc. Translations 73, 1989.

\bibitem [Qi]{Qi} S. Qian, $\varepsilon$-isometric embeddings, Proc. Amer. 
Math.  Soc. 123, 1995, 1797--1803.

\bibitem [\v Se1]{Se1} P. \v Semrl, Hyers-Ulam stability of isometries,
Houston J. Math. 24, 1998, 699--706.

\bibitem [\v Se2]{Se2} P. \v Semrl, Hyers-Ulam stability of isometries on
Banach spaces, Aequationes Math. 58, 1999, 157--162.

\bibitem [\v SV]{SV.prep} P. \v Semrl and J. V\"ais\"al\"a, Nonsurjective
nearisometries of Banach spaces, preprint.

\bibitem [Ta]{Ta} J. Tabor, Stability of surjectivity, J. Approx. Theory
105, 2000, 166--175.

\bibitem[Tr]{Tr} D.A. Trotsenko, Continuation of space quasiconformal maps
that are close to conformal maps, Sibirsk. Mat.  Zh. 28, 1989, 126--133.
(Russian)

\bibitem[Ul]{Ul} S.M. Ulam, A collection of mathematical problems,
Interscience, 1960.

\bibitem[V\"a1]{blep} J. V\"ais\"al\"a, Bilipschitz and quasisymmetric
extension properties, Ann. Acad. Sci. Fenn. Math. 11, 1986, 239--274.

\bibitem [V\"a2]{solar} J. V\"ais\"al\"a, Isometric approximation property
in euclidean spaces. Israel J. Math., to appear.

\bibitem [V\"a3]{unbdd} J. V\"ais\"al\"a, Isometric approximation property
of unbounded sets, preprint.





\end{thebibliography}
\end{document}